 \documentclass[a4paper,12pt]{article}%%book}%,landscape
 \usepackage{amsmath, amsfonts, amssymb}
 \usepackage{amscd}
 \usepackage[russian]{babel}
 \usepackage{makeidx}
 \usepackage{calrsfs}
 \usepackage{graphicx}
 \usepackage{epsfig}
 \usepackage{bbm}
\usepackage[pdftex,unicode,colorlinks,linkcolor=blue,
citecolor=red,bookmarksopen,pdfhighlight=/N]{hyperref}
 \makeindex

 \DeclareMathAlphabet{\mathbbt}{U}{bbold}{m}{n}

\begin{document}

 \makeatletter%{@}

  \newcommand{\Title}[1]{\begin{center}\large\uppercase{#1}\end{center}\par}

 \newcommand{\Author}[2]{\begin{center}\textbf{\large #1} \end{center} \medskip
                    \renewcommand{\@evenhead}{\raisebox{1mm}[\headheight][0pt]%
                    {\vbox{\hbox to\textwidth{\thepage \hfill\strut {\small #2}\hfill}\hrule}}} }

 \newcommand{\shorttitle}[1]{\renewcommand{\@oddhead}{\raisebox{1mm}[\headheight][0pt]%
                    {\vbox{\hbox to\textwidth{\strut \hfill{\small #1}\hfill\thepage}\hrule}}} }
 \headsep=2mm

  \newcommand{\l@abcd}[2]{{\begin{center}\par\medskip\normalsize\par\smallskip\noindent\hangindent5pt\hangafter=1{\bf #1}\end{center}}\par\medskip}

  \renewcommand{\section}{\@startsection{section}{1}{\parindent}
                          {3ex plus 1ex minus .2ex}{2ex plus .2ex}{\bf\centering}}
 \renewcommand{\l@section}[2]{\small\leftskip0pt\par\noindent\hangindent27pt\hangafter=1{\qquad#1} \dotfill~~#2\par}

 \renewcommand{\l@part}[2]{\normalsize\leftskip0pt\par\smallskip\noindent\hangindent17pt\hangafter=1{\bf #1} \normalsize\dotfill~~#2\smallskip\par}

 \def\subsec#1{\smallskip\textbf{#1}}

 \newcommand{\Adress}[1]{\par\bigskip\baselineskip=11pt\hangindent17pt\hangafter=0\noindent{\footnotesize#1}\par\normalsize}

 \newcommand{\UDC}[1]{\begingroup\newpage\thispagestyle{empty}\begin{flushleft}УДК #1\end{flushleft}}

 \newcommand{\Abstract}[1]{\hangindent17pt\hangafter=0\noindent{\footnotesize#1}\bigskip\par\medskip}

 \makeatother%@

 \newcommand{\affiliation}[1]{{\itshape #1}}
 \newcommand{\email}[1]{\texttt{#1}}
 %%%%%%%%%%%%%%%%%%%%%%%%%%%%%%%

 \newcommand{\bib}[2]{{\leftskip-10pt\baselineskip=11pt\footnotesize\item{}\textsl{#1}~#2\par}}
 \newcommand{\Proclaim}[1]{\smallskip{\bf#1}\sl}
 \newcommand{\proclaim}[1]{{\bf#1}\sl}
 \newcommand{\Demo}[1]{\smallskip\par{\sc #1}}
 \newcommand{\demo}[1]{{\sc #1}}
 \newcommand{\Endproc}{\rm}
 \newcommand{\Enddemo}{\rm}

%%%%%%%%%%%%%%%%%%%%%%%%%%%%%%%%%%%%%%%%%%%%%%
% Пожалуйста, не меняйте порядок команд!
%%%%%%%%%%%%%%%%%%%%%%%%%%%%%%%%%%%%%%%%%%%%%%
%
\UDC{512.5 : 517.982}

% Название статьи
\Title{Порядковые версии Теоремы Хана--Банаха и огибающие. I. Однородные функции} % обязательное поле!
%
% Авторы
\Author{Хабибуллин Б.~Н.,  Розит А.~П., Хабибуллин Ф.~Б.}
% обязательное поле!

\Abstract{Мы приводим здесь общую постановку задачи существования и построения  верхней и нижней огибающей 
для произвольной функции со значениями из пополнения упорядоченного множества ${\rm S}$ по некоторому классу функций со значениями из ${\rm S}$. Задача разбирается пока только для простейшего случая модельного класса однородных функций. Рассматриваем лишь порядково-алгебраические версии без привлечения топологии. \par}

\bigskip

\section{1. Введение. Определения и постановки задач}

Одна из классических форм  Теоремы Хана--Банаха для векторных пространств $X$ над полем вещественных чисел $\mathbb R$ гласит [1]: {\it любая  положительно однородная субаддитивная, и только такая,  функция $f\colon X\rightarrow {\mathbb R}$  равна поточечной точной верхней грани всех линейных функций $\varphi \colon X \rightarrow{\mathbb R}$, мажорируемых поточечно функцией  $f$ в том смысле, что  $\varphi (x)\leqslant f(x)$ для всех  $x\in X$}, т.\,е. функция $f$ совпадает со своей  нижней огибающей  по классу линейных функций.  
В <<Математической энциклопедии>> [2; Хана--Банаха теорема] формулировка Теоремы Хана--Банаха для векторного пространства некорректна: {\it <<В случае действительного пространства $X$ полунорму можно заменить положительно однородным функционалом, \dots>>}, т.\,е. опущено требование субаддитивности.  Тем не менее, основным ориентиром в выборе терминологии, где это возможно,  выбрана именно <<Математическая энциклопедия>> [2].
При этом, поскольку в различных источниках и у разных авторов терминология  зачастую существенно разнится,  в нашем изложении по возможности все, даже элементарные, определения, понятия и утверждения, встречавшиеся нам в литературе хотя бы раз  в различных смыслах и трактовках, приводятся полностью во избежание  разночтений. 

Дадим здесь возможную общую постановку этой проблематики, мотивированную для нас  предшествующими применениями утверждений подобного рода в теории функций [3]--[4] {\large(}см. также [5]--[8]{\large)}.

\subsec{1.1. Упорядоченные множества. Пополнение.} Пусть ${\rm S}$  ---  (частично) упорядоченное множество [1]  с отношением порядка (рефлексивным, транзитивным, антисимметричным) $\leqslant$, т.\,е. пара $({\rm S},\leqslant )$; 
$\,\geqslant\,$ и $\,>\,$ ---  соотв.\footnote{сокращение для <<соответственно>>}  {\it обратные к\/} $\,\leqslant\,$ и {\it строгому порядку\/}  $\,<\;:=\;\leqslant\, \cap \,\neq \,$. 
 
Пара $({\rm S},\leqslant)$, или множество ${\rm S}$, {\it полное снизу} (соотв. {\it сверху\/}), если для каждого непустого подмножества
 ${\rm S}_0\subset {\rm S}$ существует  {\it  точная нижняя\/} (соотв. {\it верхняя\/}) {\it граница\/} $\inf {\rm S}_0$ (соотв. $\sup {\rm S}_0$). Множество   ${\rm S}$ {\it полное,\/} 
если  ${\rm S}$ полное и снизу, и сверху.       Подмножество ${\rm S}_0\subset {\rm S}$  {\it ограничено снизу\/} (соотв. {\it сверху\/}), если 
существует элемент  $s_0 \in {\rm S}$, для которого  $s_0\leqslant s$ (соотв. $s\leqslant s_0$) для всех  $s\in {\rm S}_0$.  
Множество  ${\rm S}$  {\it порядково полное} {\it снизу\/} (соотв. {\it сверху\/})  [9]--[11]\footnote{{\it lower\/} ({\it upper\/} resp.) {\it order-complete}}, если для каждого непустого ограниченного снизу (соотв. сверху) подмножества  ${\rm S}_0\subset {\rm S}$ существует $\inf {\rm S}_0\in {\rm S}$ (соотв. $\sup {\rm S}_0\in {\rm S}$). Множество ${\rm S}$  {\it порядково полное\/,} если ${\rm S}$ порядково полное и снизу, и сверху.  

Пусть  ${\rm S}$ --- {\it порядково полное.\/} 
Если   $\inf {\rm S} $ и/или $\sup {\rm S} $ не существуют, то часто удобна и полезна операция  {\it (полу-)пополнения\/}  
порядково полного  ${\rm S} $ до полного  пут\"ем добавления  {\it  \underline{символов}\/}  $ \inf {\rm S} $ и/или $ \sup {\rm S} $, если таких элементов  первоначально  в  ${\rm S}$ нет.
Конкретнее\footnote{в  [9] использованы иные обозначения}, 
\begin{enumerate}
\item[{[$\downarrow$]}] ${\rm S}_{\downarrow} :=\{\inf {\rm S}\}\cup {\rm S}$ --- {\it полупополнение\/}, или полурасширение, множества ${\rm S}$ {\it вниз,\/} или влево;
\item[{[$\uparrow$]}] ${{\rm S}}^{\uparrow}:={\rm S}\cup \{\sup {\rm S}\}$  --- {\it полупополнение\/}, или полурасширение, множества ${\rm S}$  {\it вверх,\/} или вправо;
\item[{[$\updownarrow$]}] ${{\rm S}^{_\uparrow}_{^\downarrow}}:={\rm S}_{\downarrow}\cup {\rm S}^{\uparrow}$ --- {\it пополнение\/}, или расширение,  множества  ${\rm S}$ в порядковом смысле, 
\end{enumerate}
 где порядок  $\leqslant\,$ продолжен  естественным пут\"ем на эти пополнения, т.\,е.  $\inf {\rm S}\leqslant s\leqslant \sup {\rm S}$ {\it для всех элементов\/} 
$s\in {{\rm S}^{_\uparrow}_{^\downarrow}}$.  Очевидно, пополнения  $ {{\rm S}^{_\uparrow}_{^\downarrow}}, {{\rm S}_\downarrow}, 
{\rm S}^{\uparrow}$  с таким отношением порядка соотв. полное, полное снизу, полное сверху упорядоченные множества.
При этом, в обозначении  $\varnothing $ для пустого множества,  естественно полагать
\begin{equation*}
	\sup \varnothing :=\inf {\rm S} \quad \text{для $\varnothing\subset {\rm S}_{\downarrow}, {{\rm S}^{_\uparrow}_{^\downarrow}}$} ,
	\quad \inf \varnothing :=\sup {\rm S} \quad\text{для  $\varnothing\subset {{\rm S}^{_\uparrow}_{^\downarrow}}, 	{{\rm S}}^{\uparrow}$} .
\eqno{(0)}
\end{equation*}

\subsec{1.2. Верхняя и нижняя огибающие. Постановки задач.}
Для множеств $X,Y$ традиционно через $Y^X$ обозначаем множество всех {\it функций} 
(отображений, операторов, функционалов, форм и проч.)\footnote{Для $f\colon X\rightarrow Y$  в основном будем использовать термин {\it функция,\/} индифферентный к природе множеств $X$ и $Y$ [2; Функция].}   $f\colon X\rightarrow Y$, или $f\colon x\mapsto f(x)$, $x\in X$, или $x\mapsto f(x)$, $x\in X$, определ\"енных на $X$.  Для $X_0\subset X$ \; $f\bigm|_{X_0}$ --- {\it сужение\/ $f$ на\/} $X_0$.

Пишем  $\varphi =f$ {\it на\/} $X$,  если $\varphi (x)=f (x)$ {\it для всех\/} $x\in {X}$.  В противном случае $\varphi \neq f$ {\it на} $X$.
 Пусть $Y=\rm  {\rm S}^{_\uparrow}_{^\downarrow}$  --- пополнение  порядково полного $({\rm S}, \leqslant\,)$.  
Пишем $ \varphi \leqslant f$ {\it на\/} $X$, если $\varphi (x) \leqslant  f(x)$ {\it для всех\/ $x\in X$} и говорим, что   $\varphi$ {\it минорирует\/} $f$,
или $f$ {\it мажорирует\/} $\varphi$, {\it на\/}  $X$.  Отношение <<$f\leqslant \varphi$ {\it на\/} $X$>> определяет {\it отношение поточечного порядка на}
$\bigl({\rm S}^{_\uparrow}_{^\downarrow}\bigr)^X$. Очевидно, множество $\bigl({\rm S}^{_\uparrow}_{^\downarrow}\bigr)^X$ с отношением поточечного порядка, обозначаемого тем же символом $\,\leqslant\,$, {\it полное,\/} а именно:
{\it для произвольного   $F \subset \bigl({\rm S}^{_\uparrow}_{^\downarrow}\bigr)^X$ всегда существуют функции 
\begin{equation*}\label{siF} 
	\sup  F\colon x\mapsto \sup_{f\in  F}f(x)\in {{\rm S}^{_\uparrow}_{^\downarrow}}, \quad
	\inf F\colon x\mapsto \inf_{f\in F}f(x)\in {{\rm S}^{_\uparrow}_{^\downarrow}}, \quad x\in X,
\end{equation*}}
когда на  $X$ рассматриваются и постоянные  функции
\begin{equation*}
	\inf{\bigl({\rm S}^{_\uparrow}_{^\downarrow}\bigr)^X} 			\colon x\mapsto \inf S\in 	{\rm S}^{_\uparrow}_{^\downarrow}, \quad 
		\sup {\bigl({\rm S}^{_\uparrow}_{^\downarrow}\bigr)^X} \colon x\mapsto \sup S \in { {\rm S}^{_\uparrow}_{^\downarrow}}	, \quad x\in X,
\tag{$\star$}
	\end{equation*}
а для пустого  подмножества $\varnothing\subset \bigl({\rm S}^{_\uparrow}_{^\downarrow}\bigr)^X$  
в соответствии с соглашением  (0) определены точные  границы
\begin{equation*}\label{df:pmi}
\sup \varnothing :=\inf{\bigl({\rm S}^{_\uparrow}_{^\downarrow}\bigr)^X} \in 
\bigl({\rm S}^{_\uparrow}_{^\downarrow}\bigr)^X, \quad 
\inf \varnothing :=\sup {\bigl({\rm S}^{_\uparrow}_{^\downarrow}\bigr)^X} \in 
\bigl({\rm S}^{_\uparrow}_{^\downarrow}\bigr)^X .
\tag{$\varnothing$}
\end{equation*}

\Demo{Определение 1.}
   Пусть\/  $(S,\leqslant\,)$  порядково полное,   $f\in \bigl({\rm S}^{_\uparrow}_{^\downarrow}\bigr)^X$,  а также  $\Phi \subset \bigl({\rm S}^{_\uparrow}_{^\downarrow}\bigr)^{X}$. {\it Нижнюю\/} (соотв. {\it верхнюю\/})  {\it $\Phi$-огибающую\/}, или {\it огибающую по\/} $\Phi$  {\it  для\/}  $f$ 
{\it на\/ $X$} определяем как {\it функцию
\begin{align*}	 
{\text{\rm lE}}_{\Phi}^f&\colon	x \mapsto \sup  \bigl\{ \varphi (x) \colon  \Phi  \ni \varphi\leqslant f \text{ на }  {X}\bigr\}\in 
{{\rm S}^{_\uparrow}_{^\downarrow}}, \quad x\in {X}  
 		\\
\Bigl(\text{соотв. }
{\text{\rm uE}}^{\Phi}_f&\colon s\mapsto \inf \bigl\{ \varphi (s) \colon f \leqslant \varphi \in  \Phi   \text{ на }  X\bigr\}\in {{\rm S}^{_\uparrow}_{^\downarrow}}, \quad x\in X{\Bigr)}.
			\end{align*}
}
  	\Enddemo
		Функция $f\colon X\rightarrow Y$ с упорядоченными  $(X, \leq )$, $(Y, \leq )$ {\it возрастающая} на $X$, если для любых $x_1,x_2\in X$
из $x_1\leq  x_2$ следует $f(x_1)\leq  f(x_2)$, и {\it строго возрастающая,\/} если из $x_1< x_2$ следует $f(x_1)< f(x_2)$. Аналогично для убывания.
Функция  (строго) возрастающая или  убывающая --- ({\it строго\/}) {\it монотонная.} 	Очевидно, функции 
$f\mapsto {\text{\rm lE}}_{\Phi}^f$ и $f\mapsto {\text{\rm uE}}^{\Phi}_f$, $f\in \bigl({\rm S}^{_\uparrow}_{^\downarrow}\bigr)^X$, 
возрастающие на  полном множестве $\bigl({\rm S}^{_\uparrow}_{^\downarrow}\bigr)^X$, но, вообще говоря, не строго возрастающие.

В  приложениях роль   $X$ из Определения 1  часто играет некоторый  класс функций [4]--[6], а  ${\rm S}=\mathbb R$.
Основные  общие проблемы, диктуемые Определением 1 в ракурсе Теоремы Хана--Банаха, ---  

\Proclaim{Задача 1.}
Описать по возможности максимальный  класс функций $f\in \bigl({\rm S}^{_\uparrow}_{^\downarrow}\bigr)^X$,   равных своей нижней (верхней)\/   $\Phi$-огибающей. 
 \Endproc

\Proclaim{Задача 2.}
  Указать метод(ы)  в  той или иной мере конструктивного построения  $\Phi$-огибающих\/ ${\text{\rm lE}}_{\Phi}^f$ и\/ 
$\text{\rm uE}^{\Phi}_f$ для $f\in \bigl({\rm S}^{_\uparrow}_{^\downarrow}\bigr)^X$.
 \Endproc

В связи с приложениями к теории функций не меньший, а для определ\"енных применений  (см., например, [3]--[7]) даже  б\'ольший интерес, чем Задача 1,   представляет собой менее требовательная

\Proclaim{Задача $\bf 1^{l}$.} Описать по возможности максимальный класс 
функций $f\in \bigl({\rm S}^{_\uparrow}_{^\downarrow}\bigr)^X$, для которого\/  
${\text{\rm lE}}_{\Phi}^f  \neq   \inf{\bigl({\rm S}^{_\uparrow}_{^\downarrow}\bigr)^X}$    и/или\/ 
${\text{\rm uE}}^{\Phi}_f  \neq   \sup{\bigl({\rm S}^{_\uparrow}_{^\downarrow}\bigr)^X}$. 
\Endproc

В данной работе рассматриваем только классы $\Phi$, состоящие из однородных функций и их естественных обобщений.

Если порядково полное  ${\rm S}$ изначально дополнительно снабжено какими-либо алгебраическими операциями, согласованными с отношением порядка $\,\leqslant\,$, то продолжения этих операций на пополнения, или (полу)расширения, вообще говоря, с ограничениями, может определяется в каждом конкретном случае в зависимости от постановки рассматриваемых проблем  {\large(}см. и ср. [12; \S~4], [9; {\bf 1.3.1}]{\large)}.   

Для пары добавляемых символов $\inf {\rm S}\notin {\rm S}$ и/или $\sup {\rm S}\notin {\rm S}$ часто, в особенности, если на ${\rm S}$ используется  бинарная операция в аддитивной форме как основная, по аналогии с расширениями вещественной прямой  $\mathbb  R$ (c операцией сложения) вверх и/или вниз, используются соотв. обозначения 
$-\infty:=\inf {\rm S}$ и/или $+\infty:=\sup {\rm S}$,  а также, если на ${\rm S}$ используется  бинарная операция в мультипликативной форме как основная, по аналогии с (полу-)расширениями строго  положительного луча  $]0,+\infty[:=\{r\in \mathbb R\colon r>0\}$ (с операцией умножения) естественнее  обозначения  $0$ и/или $+\infty$.  Применяются, но не здесь, --- в зависимости от контекста и природы упорядоченного множества $(\rm S, \,\leqslant \,)$, --- также и обозначения соотв. $0, 1$; $\bf 0, 1$; $-\infty, 0$ и т.\,п.

$\mathbb N:=\{1,2, \dots\}$ и $\mathbb Z$ --- множества   {\it натуральных\/} и {\it целых чисел.\/} 
 Далее 	{\it положительность\/}  понимаем как $\, \geqslant 0$, а $\,>0\,$ --- {\it строгая положительность.\/} Аналогично  понимается 
{\it отрицательность.\/}	

Общепринятые алгебраические определения и факты --- из стандартного  университетского курса алгебры. 

\section{2. Огибающие по  однородным функциям}\label{sec:1}

\subsec{2.1. Определение и свойства однородных функций.}\label{df:hf}  Пусть $(H, \cdot )$ --- {\it полугруппа с мультипликативной формой записи,\/} т.\,е. множество с ассоциативной бинарной операцией $\boldsymbol{\cdot}\,\colon H^2\rightarrow H$, или $\boldsymbol{\cdot}\,  \colon (h_1,h_2)\mapsto h_1\cdot h_2=:h_1h_2$,    $h_1,h_2\in H$.  Множество $X$ --- {\it $H$-мно\-ж\-е\-ство,\/}  если на множестве $X$ определено   действие полугруппы $H$ (слева).  Точнее, на $H$-множестве $X$ задана операция умножения (слева) на элементы полугруппы $H$ по правилу
  $(h,x)\mapsto hx\in X$, $h\in H$, $x\in X$, с аксиомой ассоциативности
 \begin{enumerate}
	 \item[{\tt Ax0.}]\label{Ax}
		\; {\it $h_1(h_2x)=(h_1h_2)x$ для любых\/ $x\in X$ и\/ $h_1,h_2\in H$}.
	 \end{enumerate}

Если  $H$ --- полугруппа с {\it единичным элементом\/} $1$, т.\,е. {\it  моноид,\/}
 то определение $H$-множества $X$ дополняем ещ\"е одной аксиомой единичного элемента 
\begin{enumerate}
\item[\tt{Ax1}.]\label{Ax1} \; $1 x:=1\cdot x=x$ {\it для  любого элемента\/ $x\in X$}. 
	 \end{enumerate}

Пусть  $(\tt H,\cdot)$ --- ещ\"е одна, вообще говоря, другая полугруппа ${\tt H}$ с другой операцией умножения ${\boldsymbol{\cdot}\,}\colon 
{\tt (h_1,h_2) \mapsto h_1h_2\in H}$. 

\Demo{Определение 2.} Пусть $X$ --- $H$-множество, $\rm S$ --- $\tt H$-множество, а $\mathfrak h$ --- гомоморфизм полугруппы  $H$ в полугруппу $\tt H$.
Функцию $f\colon X\rightarrow {\rm S}$ называем {\it $\mathfrak h$-однородной,\/} если  
		\begin{equation*}\label{df:ph} 
		f(hx)=\mathfrak h(h)f(x) \; \text{\it для любых   $x\in X$ и $h\in H$.}
	\tag{hg}
	\end{equation*}
		Множество\/ $\mathfrak h$-однородных функций\/ $f\in {\rm S}^X$ 				обозначаем как  $\mathfrak h$-hg$\,(X)$. 
	\Enddemo
	
	\underline{Всюду далее} $H$ и ${\tt H}$ --- как минимум, {\it полугруппы,\/}    а $\mathfrak h\colon H\rightarrow {\tt H}$ --- {\it гомоморфизм полугрупп,\/}
	 $X$ --- $H$-{\it множество,\/} ${\rm S}$ --- ${\tt H}$-{\it множество,\/} а при использовании множества $\mathfrak h$-hg$\,(X)$ часто не указываем $H$-множество $X$, т.\,е. пишем просто $\mathfrak h${\rm -hg}. Кроме того, {\it всюду далее\/} в определениях и, как следствие, в утверждениях, фигурирует только образ $\mathfrak h(H)\subset {\tt H}$  полугруппы $H$. Поэтому, не умаляя общности, можем \underline{всюду считать}, что $\mathfrak h$ --- {\it эпиморфизм полугрупп.\/} При этом в случае, когда $(H,\cdot, 1)$ ещ\"е  и {\it группа,\/} то 
	гомоморфный образ группы $H$ в полугруппе ${\tt H}$ можно рассматривать как группу с ${\tt 1}:=\mathfrak h (1)$
	и $\bigl({\mathfrak h}(h)\bigr)^{-1}:=\mathfrak h (h^{-1})$, $h\in H$. Таким образом, в случае {\it группы\/} $H$, не умаляя общности, \underline{можем считать}  {\it группой} и ${\tt H}=\mathfrak h(H)$.
	
	 Далее $\rm  J$ --- множество индексов произвольной природы.  Примеры 1--7 относятся к  классическим видам однородных функций.
		
	\demo{Пример 1.}\label{exa:1} Пусть $\mathbb R_*=\mathbb R\setminus \{0\}$ --- {\it <<проколотая>> вещественная ось}, 
	$H={\tt H}=\bigl(\mathbb R_* , \boldsymbol{\cdot }, 1\bigr)$ --- мультипликативная группа c обычным умножением $\,\boldsymbol{\cdot}\,$.  
	При фиксированном $p\in \mathbb Z$  определ\"ен гомоморфизм $\mathfrak h\colon r\mapsto r^p$, $r\in \mathbb R_*$. Пусть  при каждом  $\rm j\in J$ множество 
$C_{\rm j}$ --- экземпляр вещественной оси  $\mathbb R$  и при различных $\rm j_1\neq j_2$ каждые два экземпляра   множеств $C_{\,\rm j_1}$ и $C_{\,\rm j_2}$ имеют единственную общую точку $0\in \mathbb R$. Полагаем     $X:=\cup_{\rm j\in J}\, C_{\,\rm j}$, ${\rm S}:=\mathbb R$. 
В этих соглашениях   $\mathfrak h$-однородные функции --- это  {\it однородные  степени\/ $p$  функции.\/} 
  \Enddemo

	\demo{Пример 2.}\label{exa:1+}  Пусть $\mathbb R^+:=\{r\in \mathbb R\colon r\geqslant 0\}$, $\mathbb R^+_*:=\mathbb R^+\setminus \{0\}=
	]0,+\infty[$;  $H={\tt H}=\bigl(\mathbb R^+_*, \boldsymbol{\cdot }, 1\bigr)$ --- мультипликативная группа строго положительных  
	чисел. При фиксированном $p\in \mathbb R$  определ\"ен гомоморфизм $\mathfrak h\colon r\mapsto r^p$, $r\in \mathbb R^+_*$. Пусть при каждом  $\rm j\in J$ множество 
$C_{\rm j}$ --- экземпляр   {\it положительного луча\/}  $\mathbb R^+$ и каждые    $C_{\,\rm j_1}$ и $C_{\,\rm j_2}$ 
  при  $\rm j_1\neq j_2$ имеют единственную общую точку $0\in \mathbb R^+$. Полагаем
  $X:=\cup_{\rm j\in J}\, C_{\,\rm j}$, ${\rm S}:=\mathbb R$. В этих соглашениях    $\mathfrak h$-однородные функции --- это  {\it положительно однородные степени\/ $p$ функции\/} [2; Однородная функция]. 
  \Enddemo
 
\demo{Пример 3.}\label{exa:3} $H=(\mathbb R_*,\boldsymbol{\cdot},1)$, ${\tt H}=\bigl(\mathbb R^{+}_* , \boldsymbol{\cdot},1\bigr)$ --- мультипликативные группы; при фиксированном $p\in \mathbb R$ гомоморфизм $\mathfrak h \colon r\mapsto |r|^p$, $r\in \mathbb R_*$; $X$ и ${\rm S}$ такие же, как в Примере 1.  
В этих соглашениях   $\mathfrak h$-однородные функции --- это  {\it абсолютно  однородные степени\/ $p$ функции}.
\Enddemo

\demo{Пример 4.}\label{exa:5} Здесь $r_0, {\tt r_0}\in \mathbb R_*^+$; $p\in \mathbb Z$; $H=\{r_0^n\colon n\in \mathbb Z\}$, ${\tt H}=\{ 
{\tt r_0^{\it n}}\colon n\in \mathbb Z\}$ --- мультипликативные циклические  группы;  гомоморфизм $\mathfrak h \colon r^n
\mapsto {\tt r_0^{\it np}}$,  $n\in \mathbb Z$; $X$ и ${\rm S}$ такие же, как в Примере 
2. При этом   $\mathfrak h$-однородные функции --- это  {\it ограниченно  (относительно  $H$) однородные степени\/ $p$ 
 функции\/}   [14].
\Enddemo

\demo{Пример 5.}\label{exa:6}  Пусть $n\in \mathbb N$;  $H=\bigl((\mathbb R^+_*)^n, \boldsymbol{\cdot }, (1,\cdots,1)\bigr)$ 
 --- группа с покомпонентным умножением векторов, ${\tt H}=\bigl(\mathbb R^+_*, \boldsymbol{\cdot}, 1 \bigr)$. При фиксированном $\vec{p}=(p_1,\dots, p_n)\in \mathbb R^n$  определ\"ен гомоморфизм 
$$
\vec{r}:=(r_1,\dots,r_n)\in (\mathbb R^+_*)^n=H, \quad
\mathfrak h (\vec{r}\,) :=\prod_{k=1}^n r_k^{p_k}\in \mathbb R^+_*={\tt H},
$$ 
 $X:=\mathbb R^n$ --- векторы-столбцы ${\vec x}\in X$ и $\vec{r} \cdot \vec{x}$ --- скалярное произведение, ${\rm S}:=\mathbb R$. Здесь   $\mathfrak h$-однородные функции --- это  {\it положительно однородные степени $\vec{p}$ функции $n$ переменных\/}  [2; Однородная функция]. 
  \Enddemo
	
	\demo{Пример 6.}\label{exa:4}  Пусть $\mathbb C$ --- поле  {\it комплексных чисел,\/} 
$\mathbb C_*:=\mathbb C\setminus \{0\}$ --- {\it проколотая комплексная плоскость.\/} 
 $H=(\mathbb C_*,\boldsymbol{\cdot},1)$, ${\tt H}=\bigl(\mathbb R^+_*, \boldsymbol{\cdot},1\bigr)$ --- группы; при фиксированном $p\in \mathbb R$ гомоморфизм $\mathfrak h \colon z\mapsto |z|^p$, $z\in \mathbb C_*$. 
Пусть   при   $\rm j\in J$ множество $C_{\,\rm j}$ --- экземпляр комплексной плоскости $\mathbb C$  и каждые  $C_{\,\rm j_1}$ и $C_{\,\rm j_2}$ при  $\rm j_1\neq j_2$  имеют единственную общую точку $0\in \mathbb C$; $X:=\cup_{\rm j\in J}C_{\,\rm j}$, ${\rm S}:=\mathbb C$.  При этом   $\mathfrak h$-однородные функции --- это  {\it комплексно  однородные степени\/ $p$ функции\/} [13; гл.~1, \S~3].
\Enddemo
	
	\demo{Пример 7.}\label{exa:7}  Пусть $H=\bigl((\mathbb C_*)^n,\boldsymbol{\cdot}, (1,\cdots,1)\bigr)$ с покомпонентным ум\-н\-о\-ж\-ением, ${\tt H}=\bigl(\mathbb R^+_*, \boldsymbol{\cdot},1\bigr)$, 	$\vec{p}\in \mathbb R^n$ --- мультииндекс, для $z=(z_1,\dots, z_n)\in \mathbb C_*^n$
	положим $\mathfrak h (z):= |z|^{\vec p}:=|z_1|^{p_1}\dots |z_n|^{p_n}$; 
 $X:=\mathbb C^n$, ${\rm S}:=\mathbb C$.  При этом   $\mathfrak h$-од\-н\-о\-р\-о\-д\-ные функции ---   {\it комплексно  однородные 
степени $\vec{p}$ функции $n$ комплексных переменных;\/} ср. с [13; гл.~1, \S~3].
\Enddemo

Привед\"ем пример с разными видами операций в $H$ и ${\tt H}$. 

\demo{Пример 8.}  Группы $(\mathbb R,+,0)$, $(\mathbb R_*^+, \cdot, 1)$, $p\in \mathbb R$,
 изоморфизм $\mathfrak h:=\exp^p\colon x\mapsto e^{px}$, $x\in \mathbb R$; $X=\mathbb R={\rm S}$; 
функция  $f\colon \mathbb R\rightarrow \mathbb R$ из класса $\exp^p$-hg$(\mathbb R)$, если и только если $f(h+x)=e^{ph}f(x)$ для всех $h,x\in \mathbb R$.
\Enddemo 

\underline{Другие общие примеры} можно строить для произвольной группы преобразований $H$ области определения $X$ функций $f\in {\rm S}^X$, действия-суперпозиции  $(h,f)\mapsto f\circ h$,  голоморфизма $\mathfrak h$ из $H$ в некоторую, вообще говоря, другую 
группу ${\tt H}$ преобразований  множества $X$. При этом  $f\in {\mathfrak h}\text{\rm -hg}\,(X)$, если и только если $f\circ h=f\circ 
\bigl({\mathfrak h}(h)\bigr)$ для всех $h\in H$. Частный случай --- классы функций, инвариантные относительно группы преобразований $H$, 
когда в роли  группы ${\tt H}$ ---  тривиальная  одноточечная  группа, состоящая   тождественного преобразования 
$\text{\rm id}_X$ множества $X$,  с очевидным гомоморфизмом $\mathfrak h (h)=\text{\rm id}_X$ для всех $h\in H$.   Возможности дальнейшего развития этого общего примера, Примера 8 в части различных $H$ и 
${\tt H}$ и операций на них, да  и классических Примеров 1--7 поистине неисчерпаемы.   
	
	\Demo{Определение 3.}	$\tt H$-множество $\rm S$ называем {\it упорядоченным,\/}
 если оно снабжено отношением порядка $\, \leqslant\,$ с аксиомой согласованности
\begin{enumerate}
\item[{\tt Ax2}.]\label{Ax2} {\it Для любых ${\rm s_1,s_2}\in {\rm S}$ и  $\tt h\in H$ из 
	${\rm s_1\leqslant s_2}$ следует ${\tt h}{\rm s_1}\leqslant {\tt h}{\rm s_2}$.}
\end{enumerate}
\underline{Всюду далее} ${\rm S}$ --- {\it  порядково полное  упорядоченное\/ $\tt H$-множество.\/}
		\Enddemo

Распространим  понятие $\mathfrak h$-однородной функции на $f\in \bigl({\rm S}^{_\uparrow}_{^\downarrow}\bigr)^X$. Для этого доопределим сначала действия $\tt H$ на точные грани в пополнении ${\rm S}^{_\uparrow}_{^\downarrow}$. 
Если изначально существуют $\inf {\rm S}\in {\rm S}$ и/или $\sup {\rm S} \in {\rm S}$, то  необходимости в этом нет. Иначе 
{\it для любого\/ $\tt h\in H$} полагаем
\begin{enumerate}
\item[{[$\downarrow$]}]  {\it ${\tt h}\cdot \inf {\rm S}:=\inf {\rm S}\in  {\rm S}_{\downarrow}\subset {\rm S}^{_\uparrow}_{^\downarrow}$, когда  не существует\/  $\inf {\rm S}$ в\/ ${\rm S}$;}
\item[{[$\uparrow$]}]   {\it  ${\tt h}\cdot \sup {\rm S}:=\sup  {\rm S}\in  {\rm S}^{\uparrow}\subset   {\rm S}^{_\uparrow}_{^\downarrow}$, когда  не существует\/  $\sup {\rm S}$ в\/ ${\rm S}$}. 	
\end{enumerate}
Таким образом, $\tt H$-множество ${\rm S}$ можно расширить до 
{\it $\tt H$-множества\/ $ {\rm S}^{_\uparrow}_{^\downarrow}$.} При этом  для любой функции $f\colon X\rightarrow {\rm S}^{_\uparrow}_{^\downarrow}$ 
корректно условие (hg) Определения 2, которое теперь можно дословно распространить  и на такие функции $f$. 
	При этом постоянные функции  $	\inf{\bigl({\rm S}^{_\uparrow}_{^\downarrow}\bigr)^X}$
		и 	$\sup {\bigl({\rm S}^{_\uparrow}_{^\downarrow}\bigr)^X}$, определенные в ${\rm (\star)}$,  всегда $\mathfrak h$-однородны. 
		
		\Demo{Определение 4.} Множество\/ $\mathfrak h$-однородных функций\/ $f\in {({\rm S}^{_\uparrow}_{^\downarrow})}^X$ 
		 со значениями в ${\tt H}$-множестве ${\rm S}^{_\uparrow}_{^\downarrow}$ 		в рамках  соглашений [$\downarrow$] и [$\uparrow$]  
		обозначаем через  		$\mathfrak h$-hg$_{^\uparrow}^{_\downarrow}\,(X)$. 
Далее при использовании такого множества  часто не указываем $H$-множество $X$, т.\,е. пишем просто $\mathfrak h${\rm -hg}$_{^\uparrow}^{_\downarrow}\,$.		
		\Enddemo 
		
		 Во всех Примерах  2--7  на множестве   ${\rm S}$  можно задать естественный  порядок, при котором ${\rm S}$ становится 
		порядково полным ${\mathbb R}_*^+$-множеством, а ${\rm S}^{_\uparrow}_{^\downarrow}$ --- упорядоченное  ${\mathbb R}_*^+$-множество.   Например, при ${\rm S}:=\mathbb C$ в Примерах   6--7 
		можно ввести  	<<по\-к\-о\-м\-п\-о\-н\-ентный>>  порядок:  $z_1\leqslant z_2$  означает, что     $\text{Re\,} z_1\leqslant \text{Re\,} z_2$ и $\text{Im\,} z_1\leqslant \text{Im\,} z_2$, а $\mathbb C^{_\uparrow}_{^\downarrow}$, где $\text{Re}$ и $\text{Im}$ --- операции выделения действительной и мнимой части,  получается из $\mathbb C$  добавлением двух символов
		$(-\infty)+i(-\infty):=\inf \mathbb C$ и $(+\infty)+i(+\infty):=\sup \mathbb C$; $r\bigl((\pm\infty)+i(\pm\infty)\bigr):=(\pm\infty)+i(\pm\infty)$, $r\in {\mathbb R}_*^+$.

\paragraph{{\rm 2.1.1.} Верхние и нижние границы.}\label{sss:F}  

Потребуется элементарная

\proclaim{Лемма 1.} Пусть ${\rm S_0\subset {\rm S}^{_\uparrow}_{^\downarrow}}$, ${\tt h\in H}$. 
 Тогда $\sup_{\rm s\in S_0}{\tt h} {\rm  s}\leqslant {\tt h} \sup_{\rm s\in S_0} {\rm s} $ и 
$ {\tt h} \inf_{\rm s\in S_0} {\rm s} \leqslant \inf_{\rm s\in S_0}{\tt h} {\rm  s}$. Если же\/  ${\tt H}$ --- группа, то в этих двух  неравенствах  
можно поставить знак равенства.
 \Endproc
\Demo{Доказательство.} Для ${\rm s\in S_0}$ из
${\rm s}\leqslant \sup_{\rm s\in S_0}{\rm s}$ и  аксиомы {\tt Ax2} сразу следует ${\tt h}{\rm s}\leqslant {\tt h} \sup_{\rm s\in S_0}{\rm s}$, откуда
$\sup_{\rm s\in S_0}{\tt h} {\rm  s}\leqslant {\tt h} \sup_{\rm s\in S_0} {\rm s} $. Аналогично для $\inf$.
Если ${\tt H}$ --- группа, то для $\sup$ в ${\rm S}^{_\uparrow}_{^\downarrow}$ имеем также 
\begin{multline*}
	{\tt h} \sup_{\rm s\in S_0} {\rm s}\overset{{\tt Ax1}}{=}
	{\tt h} \sup_{\rm s\in S_0} 1{\rm s}={\tt h} \sup_{\rm s\in S_0} ({\tt h^{-1}h}) {\rm s}
	\overset{{\tt Ax0}}{=}{\tt h} \sup_{\rm s\in S_0} {\tt h^{-1}}({\tt h}{\rm s})
	\\
	\bigl|\text{из доказанного выше} \bigr|\leqslant  {\tt h}\bigl({\tt h^{-1}} \sup_{\rm s\in S_0} {\tt h}{\rm s}\bigr)
	\overset{{\tt Ax0}}{=} ({\tt h}{\tt h^{-1}} )\sup_{\rm s\in S_0} {\tt h}{\rm s}=\sup_{\rm s\in S_0} {\tt h}{\rm s}.
	\end{multline*}
	Аналогично для $\inf$. Вместе с предшествующим получаем требуемое.
\Enddemo

В дополнение к Определениям  2--4 дадим

\Demo{Определение 5.} 
Функцию $f\in {({\rm S}^{_\uparrow}_{^\downarrow})}^X$ называем 
{\it $\mathfrak h$-по\-л\-у\-о\-д\-н\-о\-р\-одной снизу   (соотв. сверху) если имеем $f(hx)\leqslant \mathfrak h(h)f(x)$ 
{\rm \large(}соотв. $\mathfrak h(h)f(x) \leqslant f(hx) ${\rm \large)} для всех $x\in X$ и  $h\in H$.\/} Множество всех таких функций обозначаем 
как $\mathfrak h$-hg$_{\uparrow} (X)$ {\large(}соотв. $\mathfrak h$-hg$^{\downarrow} (X)${\large)}. 
При этом множество $X$  часто не указываем. 
Очевидно, 
$\mathfrak h\text{-hg}_{\uparrow} \cap \mathfrak h\text{-hg}^{\downarrow} =\mathfrak h$\text{-hg}$_{^\uparrow}^{_\downarrow}$.

\Enddemo

\Proclaim{Предложение 1.} Если   $F\subset  \text{$\mathfrak h${\rm -hg}$_{\uparrow}  $}$ 
{\rm \large(}соотв.\/ $F\subset  \text{$\mathfrak h${\rm -hg}$^{\downarrow}  $}${\rm \large)}, то  функция 
$\sup F\in   \text{$\mathfrak h${\rm -hg}$_{\uparrow}  $}$ {\rm \large(}соотв.\/ $\inf F\in  \text{$\mathfrak h${\rm -hg}$^{\downarrow}  $}${\rm \large)}.
В частности, если\/ $\tt H$ --- это группа, то из $F\subset  \text{$\mathfrak h${\rm -hg}$_{^\uparrow}^{_\downarrow}\, $}$ следует 
$\sup  F, \inf F\in   \text{\rm $\mathfrak h$-hg$_{^\uparrow}^{_\downarrow}\, $}$.
\Endproc

\Demo{Доказательство.} Для $F\subset  \text{$\mathfrak h${\rm -hg}$_{\uparrow}  $}$ при  $F(x):=\{f(x)\colon f\in F\}$
\begin{multline*}
	(\sup F) (hx)=\sup_{f\in F}f(hx)\leqslant \sup_{f\in F}{\mathfrak h}(h) f(x)
	=	\sup \bigl\{ {\mathfrak h}(h) {\rm s}\colon {\rm s}\in F(x) \bigr\}\\
	\Bigl|\text{Лемма 1 с ${\rm S_0}:=F(x)$}\Bigr|\leqslant 	{\mathfrak h}(h) \sup \bigl\{  {\rm s}\colon {\rm s}\in F(x) \bigr\} 
	={\mathfrak h}(h) (\sup F)(x).
\end{multline*}
 Здесь в случае $F\subset  \text{$\mathfrak h${\rm -hg}$_{^\uparrow}^{_\downarrow}\, $}$  первый знак неравенства $\,\leqslant\,$ можно заменить на равенство, а когда  ${\tt H}$ --- группа, и следующий знак  $\,\leqslant\,$ поменять на знак $\,=\,$. Аналогично для $\inf$.  Предложение 1 доказано.
\Enddemo

Для $f\in {({\rm S}^{_\uparrow}_{^\downarrow})}^X$ в обозначениях Определений 1--2, 4--5 
всегда существуют  участвующие ниже нижние   и верхние 
огибающие  
\begin{equation*}
	{\text{\rm lE}}_{\text{$\mathfrak h$-hg\,$ $}}^f\leqslant
	{\text{\rm lE}}_{\text{$\mathfrak h$\text{-hg}$_{^\uparrow}^{_\downarrow}\, $}}^f\leqslant
	{\text{\rm lE}}_{\text{$\mathfrak h${\rm -hg}$_{\uparrow}  $}}^f
		\leqslant
		f  	\leqslant 
	{\text{\rm uE}}^{\text{$\mathfrak h$\text{-hg}$^{_\downarrow}\, $}}_f
	\leqslant {\text{\rm uE}}^{\text{$\mathfrak h$\text{-hg}$_{^\uparrow}^{_\downarrow}\, $}}_f
	\leqslant 	{\text{\rm uE}}^{\text{$\mathfrak h$-hg\,$\,$}}_f \quad 
\tag{E}
\end{equation*}
на $X$, где в случае $f\in \text{$\mathfrak h$-hg\,$ \,$}$ всюду можно поставить равенства.

Из Предложения 1 и определений легко получаем 

\Proclaim{Следствие  1.}    Для $f\in {({\rm S}^{_\uparrow}_{^\downarrow})}^X$ всегда 
\begin{equation*}
{\text{\rm lE}}_{\text{\rm $\mathfrak h$-hg\,$ \,$}}^f , \, 
	{\text{\rm lE}}_{\text{$\mathfrak h$\text{\rm -hg}$_{^\uparrow}^{_\downarrow}\,$}}^f, \,
	{\text{\rm lE}}_{\text{$\mathfrak h${\rm -hg}$_{\uparrow}  \,$}}^f
	\in \text{$\mathfrak h${\rm -hg}$_{\uparrow}  \,$}
	 \quad \text{и} \quad 
	{\text{\rm uE}}^{\text{\rm $\mathfrak h$-hg\,$ \,$}}_f , \,
	{\text{\rm uE}}^{\text{$\mathfrak h$\text{\rm -hg}$_{^\uparrow}^{_\downarrow}\,$}}_f,\,
	{\text{\rm uE}}^{\text{$\mathfrak h$\text{\rm -hg}$^{_\downarrow}\,$}}_f
	\in \text{\rm $\mathfrak h$-hg$^{\downarrow}  \,$},	
\end{equation*}
а если ${\tt H}$ --- группа, то\/  
	${\text{\rm lE}}_{\text{\rm $\mathfrak h$-hg\,$ \,$}}^f , \, 
	{\text{\rm lE}}_{\text{$\mathfrak h$\text{\rm -hg}$_{^\uparrow}^{_\downarrow}\,$}}^f, \,
	{\text{\rm uE}}^{\text{\rm $\mathfrak h$-hg\,}}_f , \,
	{\text{\rm uE}}^{\text{$\mathfrak h$\text{\rm -hg}$_{^\uparrow}^{_\downarrow}\,$}}_f
	\in \text{\rm $\mathfrak h$-hg$_{^\uparrow}^{_\downarrow}\,$}$.
\Endproc

Это Следствие 1 сразу решает  Задачу 1 для ряда  классов $\Phi$: 

\Proclaim{Теорема 1.} Пусть $f\in {({\rm S}^{_\uparrow}_{^\downarrow})}^X$. Справедливы утверждения
\begin{enumerate}
	\item[{\rm [l]}] ${\text{\rm lE}}_{\text{$\mathfrak h${\rm -hg}$_{\uparrow}  $}}^f=f$ тогда и только тогда, когда 
$f	\in \text{$\mathfrak h${\rm -hg}$_{\uparrow}  $}$. 
\item[{\rm [u]}] $f={\text{\rm uE}}^{\text{$\mathfrak h$\text{\rm -hg}$^{_\downarrow}\, $}}_f$ тогда и только тогда, когда
	$f\in \text{\rm $\mathfrak h$-hg$^{\downarrow}  $}$.
\end{enumerate}
Пусть, в дополнение,\/ ${\tt H}$ --- группа. Эквивалентны три  соотношения:
{\rm 1)} $f\in \text{\rm $\mathfrak h$-hg$_{^\uparrow}^{_\downarrow}\,  $}$; \quad {\rm 2)} ${\text{\rm lE}}_{\text{$\mathfrak h$\text{\rm -hg}$_{^\uparrow}^{_\downarrow}\, $}}^f=f$; \quad
	{\rm 3)} $f={\text{\rm uE}}^{\text{$\mathfrak h$\text{\rm -hg}$_{^\uparrow}^{_\downarrow}\, $}}_f$.
	
\Endproc
Случай класса $\Phi=\text{\rm $\mathfrak h$-hg$\,  $} $ потребует дополнительных сведений.

 \demo{Замечание 1.} При доказательстве [l]--[u] дополнительные условия на полугруппы $H$ и ${\tt H}$ не накладывались. В случае же 
$\Phi = \text{\rm $\mathfrak h$-hg$_{^\uparrow}^{_\downarrow}\,  $}$ в Теореме 1  для возможности использования 
Предложения 1 и Леммы 1 предполагалось, что ${\tt H}$ --- группа.
  \Enddemo	

\paragraph{{\rm 2.1.2.} Орбиты и стационарные элементы.}   Для  элемента $x$ из $H$-мно\-ж\-е\-с\-т\-ва $X$ определяется {\it орбита\/} этого элемента 
\begin{equation}
		\text{orb}_X (x):={H}\cdot x:=:{ H}x:=\{{h}x\colon {h\in \ H}\}.
	\tag{o}
	\end{equation}
 Элемент $x\in X$ {\it стационарный,\/}  
если $\text{orb}_X(x)=\{x\}$. 

\proclaim{Основное свойство орбит.}  Для \underline{\it группы}  $H$  орбиты не зависят от выбора представителя  и  либо не пересекаются, либо совпадают. 
\Endproc

Пусть $\tt H$ --- группа. 
Тогда для любых $\rm s\in {\rm S}^{_\uparrow}_{^\downarrow}$  в рамках  соглашений [$\downarrow$] и [$\uparrow$]
определены орбиты $\text{orb}_{{\rm S}^{_\uparrow}_{^\downarrow}}({\rm s})$, а также стационарные элементы в ${\rm S}^{_\uparrow}_{^\downarrow}$. 
В частности, элементы $\inf {\rm S}\in {\rm S}^{_\uparrow}_{^\downarrow}$ при $\inf {\rm S}\notin {\rm S}$ и  
$\sup {\rm S}\in {\rm S}^{_\uparrow}_{^\downarrow}$ при $\sup {\rm S}\notin {\rm S}$ --- стационарные. Пусть  $f\in \text{\rm $\mathfrak h$-hg$_{^\uparrow}^{_\downarrow}\, $}$. 
Тривиально проверяется 
{\it
\begin{enumerate}
	\item[\rm 1)]  $f \bigl(\text{\rm orb}_X (x)\bigr)\subset \text{\rm orb}_{{\rm S}^{_\uparrow}_{^\downarrow}}f(x)$;
	в частности, если $f(x)$ --- стационарный элемент в ${\rm S}^{_\uparrow}_{^\downarrow}$, 
	то $f\bigl(\text{\rm orb}_X(x)\bigr)=\{f(x)\}$;
	\item[\rm 2)] если $\mathfrak h$ --- эпиморфизм, то  $f \bigl(\text{\rm orb}_X (x)\bigr)= \text{\rm orb}_{{\rm S}^{_\uparrow}_{^\downarrow}}f(x)$ и для  стационарного элемента $x\in X$ элемент $f(x)$ стационарный в ${\rm S}^{_\uparrow}_{^\downarrow}$.
\end{enumerate}}
 Так, для  ${\rm S}:=\mathbb R$ из Примеров 2--5  в ${\rm S}^{_\uparrow}_{^\downarrow}$ с $-\infty:=\inf \mathbb R$, $+\infty:=\sup \mathbb R$
и ${\tt H}:=\mathbb R_*^+$ --- пять орбит  $\bigl\{ \{-\infty\}, \,]-\infty,0[\,, \{0\}, \,]0,+\infty[\,, \{+\infty\} \bigr\}$ и  три стационарных элемента $\{-\infty , 0, +\infty\}$, а для  ${\rm S}:=\mathbb C$ в Примерах 6--7 орбиты  --- это либо $\{0\}$, либо  два элемента $(\pm \infty)+i(\pm\infty)$, либо 
 любой из лучей $\bigl\{re^{i\theta}\colon r\in ]0,+\infty[\bigr\}\subset \mathbb C$ 
со всевозможными $\theta \in (-\pi,\pi]$, а стационарных элементов тоже три:
 $\{0\}$ и $(\pm \infty)+i(\pm\infty)$.
	
\paragraph{{\rm 2.1.3.} Расщепление функции по орбитам.} 	Пусть $H$ --- \underline{\it  группа}. 
По Основному свойству орбит из п. 2.1.2 $H$-множество $X$ можно представить в виде объединения 
\begin{equation*}
	X=\bigcup_{\rm j\in J} \text{orb}_X( x_{\rm j}),  \quad \text{$\rm J$ --- множество идексов,}
\tag{$\rm o1$}
\end{equation*}
{\it не пересекающихся орбит\/} $\text{orb}_X( x_{\rm j})$ с выбранными по аксиоме выбора {\it представителями\/} $x_{\rm j}\in \text{orb}_X( x_{\rm j})$.
Тогда каждую орбиту $\text{orb}_X( x_{\rm j})$ можно рассматривать как $H$-множество, 
а  произвольная функция $f\in ({\rm S}^{_\uparrow}_{^\downarrow})^X$ однозначно определяется  своими сужениями 
\begin{equation*}
	f_{\rm j}:=f\bigm|_{\text{orb}_X( x_{\rm j})} \; , \quad {\rm j\in J}.
\tag{$\rm o2$}
\end{equation*}
Более того, каждая функция  $f\in \text{\rm $\mathfrak h$-hg$_{^\uparrow}^{_\downarrow}\, $}$ полностью  определяется своими значениями 
$f( x_{\rm j})$, поскольку может   быть однозначно продолжена на всю орбиту  $\text{orb}_X( x_{\rm j})$ по правилу
\begin{equation*}
	f(hx_{\rm j}):=\mathfrak h (h) f(x_{\rm j}), \quad h\in H. 
\tag{$\rm o3$}
\end{equation*}

	Следующий результат решает Задачу 1 для класса $\Phi=\text{\rm $\mathfrak h$-hg$\,  $} $. 
	
\Proclaim{Теорема 2.}	 Пусть $H$ и\/ ${\tt H}$  --- \underline{группы},  $f\in ({\rm S}^{_\uparrow}_{^\downarrow})^X$, $\text{\rm im\,}f\subset {\rm S}^{_\uparrow}_{^\downarrow}$
 --- образ функции $f$ в\/ ${\rm S}^{_\uparrow}_{^\downarrow}$. Справедливы следующие три утверждения. 
\begin{enumerate}
	\item[{\rm [l]}] ${\text{\rm lE}}_{\text{$\mathfrak h${\rm -hg}$\,  $}}^f=f$ тогда и только тогда, когда  выполнено  одно из следующих взаимоисключающих двух условий:
		\begin{enumerate}
		\item[{\rm [l1]}] $\text{\rm im\,}f\subset {\rm S}^{\uparrow}$ и\/ $f\in \text{\rm $\mathfrak h$-hg$_{^\uparrow}^{_\downarrow}\, $}$,
		\item[{\rm [l2]}] $\inf {\rm S}\notin {\rm S}$ и\/  $f=	\inf{\bigl({\rm S}^{_\uparrow}_{^\downarrow}\bigr)^X}$ на $X$\/
		{\rm {\rm \large(}см. ($\star$) в 1.2{\rm \large)}}.
			\end{enumerate}
\item[{\rm [u]}] $f={\text{\rm uE}}^{\text{$\mathfrak h${\rm -hg}$\,  $}}_f$ тогда и только тогда, когда  выполнено одно из следующих взаимоисключающих двух условий:
\begin{enumerate}
		\item[{\rm [u1]}] $\text{\rm im\,}f\subset {\rm S}_{\downarrow}$ и\/ 
		$f\in \text{\rm $\mathfrak h$-hg$_{^\uparrow}^{_\downarrow}\, $}$,
		\item[{\rm [u2]}] $\sup  {\rm S}\notin {\rm S}$ и $f=	\sup{\bigl({\rm S}^{_\uparrow}_{^\downarrow}\bigr)^X}$ на $X$\/
			{\rm {\rm \large(}см. ($\star$) в 1.2{\rm \large)}}.
			\end{enumerate}
	\item[{\rm [lu]}] ${\text{\rm lE}}_{\text{$\mathfrak h${\rm -hg}$\,  $}}^f=f={\text{\rm uE}}^{\text{$\mathfrak h${\rm -hg}$\,  $}}_f$, если и только если $f\in \text{\rm $\mathfrak h$-hg$\, $}$.
\end{enumerate}
	\Endproc
	
	\Demo{Доказательство.} Утверждение [lu]=[l]$\cap$[u] следует из [l] и [u]. Утверждение [u] доказывается так же, как и [l], которое и докажем.
	Сначала необходимость. Если ${\text{\rm lE}}_{\text{$\mathfrak h${\rm -hg}$\,  $}}^f=f$,  $\inf {\rm S}\notin {\rm S}$, а $f$ принимает значение $\inf {\rm S}$ хотя бы раз, то множество функций $\varphi \in \text{$\mathfrak h${\rm -hg}$\,  $}$, мажорируемых функцией $f$, пусто. Следовательно, по Определению 1  имеем $f=	\inf{\bigl({\rm S}^{_\uparrow}_{^\downarrow}\bigr)^X}$ 	{\rm {\rm \large(}см. ($\varnothing$) в 1.2{\rm \large)}}, т.\,е. [l2]. 
	Пусть теперь $\inf {\rm S}\in {\rm S}$ или функция $f\in ({\rm S}^{_\uparrow}_{^\downarrow})^X$ нигде не принимает значение $\inf {\rm S}$, что означает выполнение условия $\text{\rm im\,}f\subset {\rm S}^{\uparrow}$.   Тогда  
	по Следствию 1 в части, когда ${\tt H}$ --- группа, имеем $f={\text{\rm lE}}_{\text{\rm $\mathfrak h$-hg\,}}^f 	\in \text{\rm $\mathfrak h$-hg$_{^\uparrow}^{_\downarrow}\,  $}$, 
	т.\,е. [l1].

	Теперь достаточность. Если выполнено [l2], то множество функций $\varphi \in \text{$\mathfrak h${\rm -hg}$\,  $}$, мажорируемых функцией 
	$f=	\inf{\bigl({\rm S}^{_\uparrow}_{^\downarrow}\bigr)^X}$, пусто. Отсюда ввиду ($\varnothing$) из  1.2 нижняя огибающая  
	${\text{\rm lE}}_{\text{$\mathfrak h${\rm -hg}$\,  $}}^f$ всюда на $X$ принимает значение $\inf {\rm S}$, т.\,е. в этом случае 
	${\text{\rm lE}}_{\text{$\mathfrak h${\rm -hg}$\,  $}}^f=\inf{\bigl({\rm S}^{_\uparrow}_{^\downarrow}\bigr)^X}$.
	Допустим теперь, что  выполнено [l1]. Пусть $x_0\in X$ и  ${\rm s}_0 \in {\rm S}$ --- произвольный элемент, для которого ${\rm s}_0\leqslant f(x_0)$. 
	Будем пользоваться расщеплением функции $f$ по непересекающимся орбитам, описанном в начале п.~2.1.3 посредством $\rm (o1)$--$\rm (o2)$ . Здесь мы уже пользуемся тем, что $H$ --- \underline{\it группа}.    Найд\"ется индекс ${\rm j}_0$, для которого   $x_0\in \text{orb}_X(x_{{\rm j}_0})$. Поскольку в этом случае   	$\text{orb}_X(x_0) =\text{orb}_X(x_{{\rm j}_0})$, то, используя переиндексацию,  можем считать, что $x_{{\rm j}_0}=x_0$. Для  каждой из остальных  орбит  $\text{orb}_X(x_{{\rm j}})$ с ${\rm j}\neq {\rm j}_0$ выберем, вновь используя аксиому выбора, произвольный элемент ${\rm s}_{\rm j}\in {\rm S}$, удовлетворяющий 	условию ${\rm s}_{\rm j} \leqslant f(x_{\rm j})$. Построение функции $\varphi \in \text{$\mathfrak h${\rm -hg}$\,  $}$, мажорируемой  функцией 	$f=	\inf{\bigl({\rm S}^{_\uparrow}_{^\downarrow}\bigr)^X}$, а также удовлетворяющей условию $\varphi (x_0)={\rm s}_0\leqslant 
	f(x_0)$, провед\"ем по схеме ${\rm (o3)}$, а именно: положим $\varphi (x_{\rm j})={\rm s}_{\rm j}$ для всех ${\rm j\in J}$ и продолжим $\varphi$ с элемента $x_{\rm j}$ на всю орбиту $\text{orb}_X(x_{\rm j})$ по правилу  $\varphi (hx_{\rm j}):=\mathfrak h (h)\varphi (x_{\rm j})=\mathfrak h (h){\rm s}_{\rm j}$, $h\in H$. Построенная функция  $\varphi \in \text{$\mathfrak h${\rm -hg}$\,  $}$ минорирует функцию $f$, откуда, в силу произвола в выборе ${\rm s}_0\in {\rm S}$ с  ${\rm s}_0\leq f(x_0)$, получаем  
		$\sup \bigl\{ \varphi (x_0) \colon  \varphi \in \text{$\mathfrak h${\rm -hg}$\,  $}, \, \varphi  \leqslant f  \text{ на $X$} \bigr\}=f(x_0)$.
		По Определению 1 это, ввиду произвола в выборе элемента $x\in X$,  означает ${\text{\rm lE}}_{\text{$\mathfrak h${\rm -hg}$\,  $}}^f=f$, что и требовалось.
	\Enddemo

\demo{Замечание 2.} На {\it полугруппы\/} $H$\/ и ${\tt H}$  при доказательстве необходимости накладывалось только одно дополнительное 
  условие:   ${\tt H}$ --- {\it группа}.\/ Напротив, при доказательстве достаточности применялось тоже только одно 
	дополнительное    условие, но другое:    ${H}$ --- {\it группа}, из которого, впрочем, следует, что ${\tt H}=\mathfrak h (H)$ --- группа. 
  \Enddemo	

\paragraph{{\rm 2.1.4.} Регуляризованные миноранта и  мажоранта.}
Наряду с (E) из п.~2.1.1 дадим ещ\"е один  способ конструирования миноранты и мажоранты функции   $f\in {({\rm S}^{_\uparrow}_{^\downarrow})}^X$. 
Для этого, в предположении, что ${\tt H}$ --- это \underline{\textit{группа}}, определим возрастающие на ${({\rm S}^{_\uparrow}_{^\downarrow})}^X$ функции $f\mapsto f_{\wedge}$ и  $f\mapsto f^\vee$, \, $f\in {({\rm S}^{_\uparrow}_{^\downarrow})}^X$, действующие по правилу
 \begin{equation*}
	 f_{\wedge} (x):= \inf_{h\in H}  \,\bigl(\mathfrak h (h)\bigr)^{-1} f(hx) , \quad 
	f^{\vee} (x):= \sup_{h\in H}  \,\bigl(\mathfrak h (h)\bigr)^{-1} f(hx) , \quad x\in X,
\end{equation*}
--- соотв. {\it регуляризованная миноранта\/} и {\it регуляризованная  мажоранта\/} функции $f$.
 Для {\it моноида\/} $H$, очевидно,  $f_{\wedge}\leqslant f\leqslant f^{\vee}$ на $X$.

\demo{Замечание 3.} Любую полугруппу $H$ без единичного элемента можно превратить в моноид, просто присоединив формальный единичный элемент $1$ и определив $1h := h =: h1$ для всех $h \in H$. При этом, очевидно,  гомоморфизм $\mathfrak h$ корректно продолжается на моноид $H\cup \{1\}$  по правилу ${\mathfrak h}(1)={\tt 1 \in H}$ --- здесь группа.
  \Enddemo

\Proclaim{Предложение 3.} Всегда $f_{\wedge}\in \text{$\mathfrak h${\rm -hg}$^{\downarrow}  $}$ и $f^{\vee} \in 
\text{\rm $\mathfrak h$-hg$_{\uparrow}  $}$. 
Если же, в дополнение,\/  $H$ --- группа, то  $f_{\wedge}, \,f^{\vee} \in 
\text{$\mathfrak h${\rm -hg}$_{^\uparrow}^{_\downarrow}\, $}$.
\Endproc
\Demo{Доказательство.} Для произвольного $h_0\in H$ имеем
\begin{multline*}
	f_{\wedge} (h_0x)=\inf_{h\in H}\bigl(\mathfrak h (h)\bigr)^{-1} f\bigl(h(h_0x)\bigr) 
	=\inf_{h\in H} \mathfrak h (h_0) \bigl(\mathfrak h (h)\mathfrak h (h_0)\bigr)^{-1} f\bigl((hh_0)x\bigr) \\
	=\inf_{h\in H} \mathfrak h (h_0) \bigl(\mathfrak h (hh_0)\bigr)^{-1} f\bigl((hh_0)x\bigr)
	\\
\Bigl|\text{ Лемма 1 при $S_0:=\bigl\{\mathfrak h (hh_0)\bigr)^{-1} f\bigl((hh_0)x\bigr)\colon h\in H\bigr\}$ с $\mathfrak h (h_0)$ вместо ${\tt h}\,$}\Bigr| \\
= \mathfrak h (h_0)  \inf_{h\in H}  \bigl(\mathfrak h (hh_0)\bigr)^{-1} f\bigl((hh_0)x\bigr)
=\mathfrak h (h_0)  \inf_{h\in Hh_0}  \bigl(\mathfrak h (h)\bigr)^{-1} f(hx)\\
\geqslant  \mathfrak h (h_0)  \inf_{h\in H}  \bigl(\mathfrak h (h)\bigr)^{-1} f\bigl((h)x\bigr)= \mathfrak h (h_0) f_{\wedge} (x)
\end{multline*}
ввиду $Hh_0:=\{hh_0\colon h\in H\}\subset H$ для всех $x\in X$. Отсюда сразу следует $f_{\wedge}\in \text{$\mathfrak h${\rm -hg}$^{\downarrow}  $}$. 
Если же $H$ --- \underline{\it группа,\/}  то $Hh_0=H$ и знак неравенства $\,\geqslant\,$ здесь можно заменить на $\,=\,$.
Аналогично рассматривается и $f^{\vee}$.   Предложение 3 доказано. 
\Enddemo

\demo{Замечание 3.}  Из доказательства нетрудно видеть, что для выполнения $f_{\wedge}\in 
\text{$\mathfrak h$-hg$_{^\uparrow}^{_\downarrow}\, $}$ достаточно требовать, чтобы в $H$ существовали правый единичный элемент и правый обратный  для любого элемента из  $H$. Для $f^{\vee} \in \text{$\mathfrak h$-hg$_{^\uparrow}^{_\downarrow}\, $}$ то же самое, но  с левыми.
\Enddemo

Следующий результат решает Задачу 2  для различных классов $\Phi$ однородных функций.

\Proclaim{Теорема 3.} Пусть  $H,\, {\tt H} $  --- группы,  $f\in {({\rm S}^{_\uparrow}_{^\downarrow})}^X$.  Тогда 
\begin{enumerate}
	\item [{\rm I.}] 
	${\text{\rm lE}}_{\text{$\mathfrak h${\rm -hg}$_{^\uparrow}^{_\downarrow}\, $}}^f
	=f_{\wedge}\leqslant f\leqslant f^{\vee}= 
	{\text{\rm lE}}^{\text{$\mathfrak h${\rm -hg}$_{^\uparrow}^{_\downarrow}\, $}}_f$ на $X$.
			
	\item[{\rm II.}]  ${\text{\rm lE}}_{\text{$\mathfrak h${\rm -hg}$\, $}}^f 	=f_{\wedge}$ на\/ $X$ тогда и только тогда, когда 
 выполнено одно из двух взаимоисключающих условий: 
\begin{itemize}

\item[{\rm II1.}] $\text{\rm im}\, f_{\wedge}\subset {\rm S}^{\uparrow}$;

\item[{\rm II2.}]   $\inf {\rm S}\notin {\rm S}$ и\/  $f_{\wedge}=	\inf{\bigl({\rm S}^{_\uparrow}_{^\downarrow}\bigr)^X}$ 
		на $X$\/ {\rm {\rm \large(}см. ($\star$) в 1.2{\rm \large)}}.
				 \end{itemize}

	\item[{\rm III.}]  $f^{\vee}= {\text{\rm lE}}^{\text{$\mathfrak h${\rm -hg}$\, $}}_f$ на $X$   тогда и только тогда, когда 
 выполнено одно из двух взаимоисключающих условий: 
\begin{itemize}
\item[{\rm III1.}] $\text{\rm im}\, f^{\vee} \subset {\rm S}_{\downarrow}$;
\item[{\rm III2.}]   $\sup {\rm S}\notin {\rm S}$ и\/  $f^{\vee}=	\sup{\bigl({\rm S}^{_\uparrow}_{^\downarrow}\bigr)^X}$ на $X$\/
		{\rm {\rm \large(}см. ($\star$) в 1.2{\rm \large)}}.
				 \end{itemize}
	\end{enumerate}
\Endproc
\Demo{Доказательство.} I. Рассмотрим регуляризованную  миноранту $f_{\wedge} \leqslant f$. По Предложению  3  и Определению 1 
$f_{\wedge} \leqslant {\text{\rm lE}}_{\text{$\mathfrak h${\rm -hg}$_{^\uparrow}^{_\downarrow}\, $}}^f$ на $X$. Кроме того, для произвольной функции
$\varphi \in {\text{$\mathfrak h${\rm -hg}$_{^\uparrow}^{_\downarrow}\, $}}$, мажорируемой функцией $f$,  при всех  $x\in X$ по Определению 4 имеем 
\begin{equation*}
\varphi (x)=\bigl(\mathfrak h (h)\bigr)^{-1} \varphi (hx) \overset{{\tt Ax2}}{\leqslant} 
\bigl(\mathfrak h (h)\bigr)^{-1} f (hx)  \quad\text{для всех  $h\in H$}.
\end{equation*}
Отсюда 
	$\varphi (x)\leqslant \inf_{h\in H}\bigl(\mathfrak h (h)\bigr)^{-1} f (hx) = f_{\wedge} (x) \quad \text{для всех $x\in X$}$.
Это означает, что ${\text{\rm lE}}_{\text{$\mathfrak h${\rm -hg}$_{^\uparrow}^{_\downarrow}\, $}}^f \leqslant  f_{\wedge}$ на $X$, и доказывает I
для $f_{\wedge}$.

Аналогично устанавливается равенство для $f^{\vee}$ из I.

II. {\it Необходимость.\/} Если $\inf {\rm S}\notin {\rm S}$ и $f_{\wedge}(x_0)=\inf {\rm S}$, то, ввиду   доказанного в предыдущем п.~I,  
${\text{\rm lE}}_{\text{$\mathfrak h${\rm -hg}$\, $}}^f 	\leqslant f_{\wedge}$ на $X$, т.\,е. 
${\text{\rm lE}}_{\text{$\mathfrak h${\rm -hg}$\, $}}^f (x_0)=\inf {\rm S}$.  Отсюда множество функций $\varphi \in \text{$\mathfrak h${\rm -hg}$\, $}$, мажорируемых функцией $f$, а значит и $f_{\wedge}$, пусто и $\inf{\bigl({\rm S}^{_\uparrow}_{^\downarrow}\bigr)^X}={\text{\rm lE}}_{\text{$\mathfrak h${\rm -hg}$\, $}}^f 	=f_{\wedge}$. В противном случае оста\"ется только  ситуация II2, т.\,е. $\text{\rm im}\, f_{\wedge}\subset {\rm S}^{\uparrow}$. 

{\it Достаточность.\/} Если выполнено II2, то из п.~I сразу следует ${\text{\rm lE}}_{\text{$\mathfrak h${\rm -hg}$\, $}}^f\leqslant 
f_{\wedge}=\inf{\bigl({\rm S}^{_\uparrow}_{^\downarrow}\bigr)^X}$ на $X$, что и нужно. 

Пусть теперь   $\text{im\,}f_{\wedge} \subset {\rm S}^{\uparrow}$. По  Предложению 3 \,  $f_{\wedge} \in 
\text{$\mathfrak h${\rm -hg}$_{^\uparrow}^{_\downarrow}\, $}$, т.\,е. выполнено условие [l1] Теоремы 2 для функции 
$f_{\wedge}$ вместо $f$. Следовательно, по  Теореме 2 имеем ${\text{\rm lE}}_{\text{$\mathfrak h${\rm -hg}$\, $}}^f 	=f_{\wedge}$ на\/ $X$.  

III доказывается аналогично. 
\Enddemo

Отметим далее некоторые простейшие  факты об однородных функциях, которые могут быть полезными.
  
\paragraph{{\rm 2.1.5.} Умножение  функций слева.}
Для произвольной функции $f\in {({\rm S}^{_\uparrow}_{^\downarrow})}^X$ определено умножение (слева) функции $f$ на элементы ${\tt h\in H}$, а именно:
${\tt h}f\colon x\mapsto  {\tt h}f (x)$, $x\in X$.

\Proclaim{Предложение 4.} Если полугруппа ${\tt H}$  коммутативна, то для любого элемента ${\tt h\in H}$ из\/
   $f\in  \text{$\mathfrak h${\rm -hg}$^{\uparrow}  $}$  {\rm \large(}соотв.\/ $f\in   \text{$\mathfrak h${\rm -hg}$_{\downarrow}  $}${\rm \large)}
	следует  	${\tt h}f\in  \text{$\mathfrak h${\rm -hg}$^{\uparrow}  $}$   	{\rm \large(}соотв.\/ ${\tt h}f\in   \text{$\mathfrak h${\rm -hg}$_{\downarrow}  $}${\rm \large)}.  В частности, для    $f\in {\mathfrak h}$-$\text{\rm hg\,}  $ получаем\/  ${\tt h}f\in {\mathfrak h}$-$\text{\rm hg\,}  $.
\Endproc
\Demo{Доказательство.} 
$\bigl(f\in  \text{$\mathfrak h${\rm -hg}$^{\uparrow}  $}\bigr)$ 
$\Longrightarrow$ 
$({\tt h}f)(hx)\leqslant {\tt h} {\mathfrak h}(h) f(x) =\bigl|\text{коммутативность}\bigr|= {\mathfrak h}(h) {\tt h} f(x)={\mathfrak h}(h) ({\tt h} f)(x)$,   
$x\in X$, $h\in H$. Аналогично для 
$f\in   \text{$\mathfrak h${\rm -hg}$_{\downarrow}  $}$. Отсюда заключение  и для $f\in {\mathfrak h}$-$\text{\rm hg\,}  $. 
\Enddemo

	\paragraph{{\rm 2.1.6.} Локализуемость (полу)однородных функций.}  Покажем, что при  определ\"енных ограничениях на (полу)группу 
	$H$  условия (полу)однородности из Определений 2 и 4--5  
		можно проверять лишь для  порождающей части\/ $H$. Напомним, что подмножество 
$H'\subset H$ {\it порождает полугруппу\/} $H$, если множество всех конечных произведений вида $h_1\cdots h_n$ с $h_k\in H'$  
при каждом $k=1,\dots, n$  совпадает с $H$, и   {\it порождает группу} $H$, если множество всех конечных произведений того же вида, но с 
  $h_k\in H'$ или $h_k^{-1}\in H'$ при каждом $k=1,\dots, n$, --- это в точности группа $H$.

\Proclaim{Теорема 4.} Пусть $f\in {({\rm S}^{_\uparrow}_{^\downarrow})}^X$. Если  для любой неодноточечной  орбиты $\text{\rm orb}_X(\cdot) $ 
из\/ {\rm (o) п.~2.1.2}  найд\"ется подмножество  $H' \subset H$, порождающее \underline{полугруппу} $H$, с которым 
\begin{equation*}
	f(hx)\leqslant \mathfrak h (h) f(x)
\tag{\rm sh}
\end{equation*}
 для всех $h\in H'$ и $x\in \text{\rm orb}_X(\cdot)$,   то\/ {\rm (sh)} выполнено для всех $h\in H$ и $x\in X$. 
Аналогично при замене $\,\leqslant\,$ в\/ {\rm (sh)}  на $\,\geqslant\,$ или $\,=\,$.

Если   для любой неодноточечной  орбиты $\text{\rm orb}_X(\cdot) $ найд\"ется подмножество  $H' \subset H$, порождающее \underline{группу} $H$, с которым $f(hx)= \mathfrak h (h) f(x)$ для всех $h\in H'$ и $x\in \text{\rm orb}_X(\cdot)$, то
$f\in \text{$\mathfrak h${\rm -hg}$_{^\uparrow}^{_\downarrow} (X)$}$.
\Endproc

\Demo{Доказательство.} Пусть сначала $H$ --- полугруппа. Тогда в  условиях Теоремы 4 
для произвольного $h\in H$ найдутся $h_1, \dots , h_n\in H'$,  $n\in \mathbb N$, с которыми $h=h_1\cdots h_n$. Отсюда $f (hx)=
f(h_1\cdots h_nx)\overset{{\tt Ax2}, {\rm (sh)}}{\leqslant} \mathfrak h (h_1)\cdot f(h_2\cdots h_nx) \overset{{\tt Ax2}, {\rm (sh)}}{\leqslant} \cdots \overset{{\tt Ax2}, {\rm (sh)}}{\leqslant}  \mathfrak h (h_1)\cdots \mathfrak h(h_n)f(x)= \mathfrak h (h_1\cdots h_n)f(x)=\mathfrak h (h)f(x)$ для любых $h\in H$  в произвольной орбите   $\text{\rm orb}_X(\cdot)\ni x$, а значит и любых   $x\in X$, что и требуется. Для $\,\geqslant\,$ и $\,=\,$ аналогично.

Пусть теперь $H$ --- группа, а значит и ${\tt H}=\mathfrak h (H)$ --- группа.
 Для фиксированной орбиты $\text{\rm orb}\,(\cdot)$ с соответствующим порождающим $H$ множеством $H'$ 
рассмотрим элемент  $h\in H$, для которого  $h^{-1}\in H'$. Тогда для $x\in \text{\rm orb}\,(\cdot)$ 
имеем  $f(x)= f(h^{-1}hx) = \mathfrak h(h^{-1}) f(hx)=\bigl(\mathfrak h (h)\bigr)^{-1}f(hx)$, т.\,е. $f(hx)=\bigl(\mathfrak h (h)\bigr)f(x)$ 
для любых   $x\in \text{\rm orb}_X(\cdot)$ и $h^{-1}\in H'$. Таким образом, выполнено (sh) с равенством $\,=\,$ вместо неравенства $\,\leqslant\,$
для всех $h$ из множества $H'\cup (H')^{-1}\subset H$, \textit{порождающего $H$ как полугруппу.\/} 
Оста\"ется воспользоваться доказанным. 

 \Demo{Пример 9.} При $H=\mathbb R_*^+$ с обычной структурой мультипликативной группы  достаточно проверять 
$\mathfrak h$-однородность на сколь угодно коротких непустых открытых интервалах 
$\bigl]r_1,r_2\bigr[\subset \mathbb R_*^+$, поскольку, как легко показать, любой такой интервал порождает группу  $\mathbb R_*^+$.

\demo{Заключительный комментарий.} Содержательные утверждения по Задаче $1^{\rm l}$ легко следуют из Теорем 1--3. Другие варианты решения этой Задачи, не опирающиеся на эти Теоремы нам пока неизвестны. Оста\"ется неясным также вопрос --- в какой мере установленные здесь результаты для  группы  $H$ могут быть перенесены на полугруппы $H$.\/ Скорее всего всегда есть контрпримеры. Дальнейшие перспективы --- рассмотрение $\Phi$-огибающих для других классов $\Phi$, как-то: (полу=суб- или супер-)аддитивных, с условиями типа выпуклости--вогнутости,  связанных с операциями
 $\sup$ или $\inf$ (идемпотентные, или тропические, версии) и проч. Возможны вариации и в рамках структуры множества $X$ --- теоретико-множественной, алгебраической, топологической, геометрической, порядковой и далее.  Объ\"ем таких  исследований необозрим. При этом представляется важным  изучение подобных вопросов для функций со значениями именно в пополнении (расширении) ${\rm S}^{_\uparrow}_{^\downarrow}$ в связи с приложениями {\large(}см., например, [3]--[6] --- в теории функций, [15]  --- в оптимизации{\large)}.
\Enddemo

% Список литературы:
 \section{Литература}

 \begin{enumerate}
 \itemsep=0pt
 \parskip=0pt

 \bib{Акилов Г.\,П.~Кутателадзе С.\,С.}{Упорядоченные векторные пространства.---Новосибирск: Наука, 1978.}
 
\bib{}{Математическая энциклопедия. М.: <<Советская энциклопедия>>, 1977--1985.}

\bib{Хабибуллин~Б.\,Н.}{Двойственное представление суперлинейных функционалов и его применения в теории функций. I~/\!/ Изв. РАН, сер. матем.---2001.---Т.~65.---\textnumero~4---С.~205--224.}

\bib{Хабибуллин~Б.\,Н.}{Двойственное представление суперлинейных функционалов и его применения в теории функций. II~/\!/ Изв. РАН, сер. матем.---2001.---Т.~65.---\textnumero~5---С.~167--190.}

\bib{Хабибуллин~Б.\,Н.}{Двойственное представление суперлинейных функционалов~/\!/ В сб. статей <<Комплексный анализ, дифференциальные уравнения,численные методы и приложения.Часть I. Комплексный анализ.>> Уфа.---1996.---УНЦ РАН. Институт математики с ВЦ.---C.~122--131.}

\bib{Хабибуллин~Б.\,Н.}{Применения в комплексном анализе двойственного представления функционалов на векторных решетках~/\!/ Математический форум. Т.~4. Исследования по математическому анализу, дифференциальным уравнениям и 
их приложениям.---Владикавказ: ЮМИ ВНЦ РАН и РСО-А.---2010. С.~102--116 (Итоги науки. Юг России).} 

\bib{Картак~В.\,В., Хабибуллин~Б.\,Н.}{Двойственное представление функционалов на проективных пределах векторных решеток~/\!/ Теория функций, ее приложения и смежные вопросы. Труды Математического центра имени Н. И. Лобачевского. Материалы Девятой международной Казанской летней научной школы-конференции  (Казань, 1--7 июля 2009 г.). Казанское математическое общество.---Т. 38.---2009.---С. 146--148.}

\bib{Хабибуллин~Б.\,Н.}{Аналоги теоремы Хана-Банаха для (полу)групп: построение нижней огибающей~/\!/  Материалы Международной конференции <<Алгебра и математическая логика: теория и приложения>>. Казань, 2--6 июня 2014 г. Казанский (Приволжский) федеральный университет.---2014.---С. 75--76.}

\bib{Кусраев А.\,Г., Кутателадзе С.\,С.}{Субдифференциальное исчисление. Теория и приложения.---М.: Наука, 2007.}

\bib{{\rm S}imons~{\rm S}.}{From Hahn-Banach to Monotonicity.---Berlin: {\rm S}pringer {\rm S}cience+Bu\-s\-i\-n\-ess Media B.V. Lect. Notes in Math.---V.~1963, 2008.} 

\bib{Borwein~J.\, M. and Vanderwerff~J.\,D.}{Convex Functions: Constructions,
Cha\-r\-a\-c\-t\-e\-r\-i\-zations and Counterexamples.--- N.~Y.: Cambridge University Press, 2010.}

\bib{Лейхтвейс~К.}{Выпуклые множества. М.: Наука, 1985.}

\bib{Лелон~Л., Груман~Л.}{Целые функции многих комплексных переменных. М.: Мир, 1989.}

\bib{Schlude K.}{Bemerkung zu beschr\"ankt homogenen Funktionen~/\!/   
Elemente der Mathematik.---1999.---V.~54.---P.~30--31.}

\bib{Dinha~N., Ernst~ E.,  L\'opezc~M.\,A.,  Volled~M.}{An approximate Hahn--Banach theorem for positively homogeneous functions~/\!/   
Optimization: A Journal of Math. Programming and Operations Research---2015.---V.~64, 5.---P.~1321--1328.}

 \end{enumerate}

 \Adress{\textsc{Хабибуллин Булат Нурмиевич}\\
 Башкирский государственный университет\\
 заведующий кафедрой высшей алгебры и геометрии, профессор\\
 РФ, 450076, г.\,Уфа, ул. З.~Валиди, 32\\
 E-mail: \verb"Khabib-Bulat@mail.ru"}

\Adress{\textsc{Розит Алексей Петрович}\\
 аспирант кафедры высшей алгебры и геометрии\\
  E-mail: \verb"Rozit@mail.ru"}

\Adress{\textsc{Хабибуллин Фархат Булатович}\\
  доцент кафедры высшей алгебры и геометрии\\
  E-mail: \verb"KhabibullinFB@list.ru"}

 \bigskip
 \begin{center}
 ORDER VERSIONS OF THE HAHN--BANACH THEOREM\\ AND  ENVELOPES. I. HOMOGENEOUS FUNCTIONS
 \end{center}
 \begin{center}
  Bulat~N.~Khabibullin, Aleksey~P.~Rozit, Farkhat~B.~Khabibullin
	\end{center}
\Abstract{We present here a general formulation of the problem of existence and co\-n\-s\-t\-r\-u\-c\-t\-ion of upper and lower envelope for a  function with values in a completion of or\-d\-e\-r\-ed set {\rm S} for a certain class of functions with values in {\rm S}. The task is parsed on\-ly for the simplest case of model class of homogeneous functions. We consider only order-algebraic versions without the involvement of the topology.\par}

\end{document}